\documentclass[brochure,english,12pt]{bourbaki} 
\usepackage[matrix,arrow,curve]{xy}
\usepackage{amssymb,amsfonts,amsmath,footnote}
\def\CC{{{\mathbb{C}}}}
\def\AA{{{\mathbb{A}}}}
\def\FF{{{\mathbb{F}}}}
\def\GG{{{\mathbb{G}}}}
\def\PP{{{\mathbb{P}}}}
\def\QQ{{{\mathbb{Q}}}}
\def\ZZ{{{\mathbb{Z}}}}
\def\QQ{{{\mathbb{Q}}}}

\def\cF{{{\mathcal{F}}}}
\def\cE{{{\mathcal{E}}}}
\def\cM{{{\mathcal{M}}}}
\def\cI{{{\mathcal{I}}}}
\def\cL{{{\mathcal{L}}}}
\def\cO{{{\mathcal{O}}}}
\def\ss{{{\mathrm{ss}}}}
\def\fg{{\mathrm{fg}}}
\def\modfg{{{\cM_\fg}}}
\def\ord{{{\mathrm{ord}}}}

\def\hatord{{{\hat{\cM}_\ord}}}
\def\mss{{{\cM_\ss}}}
\def\hatss{{{\hat{\cM}_\ss}}}
\def\Ell{{{e\ell\ell}}}
\def\mell{{{\bar{\cM}_\Ell}}}

\def\CP{{{\CC\mathrm{P}}}}
\def\longr{{\longrightarrow\ }}
\def\tmf{{{\mathbf{tmf}}}}

\def\mstring{{{\mathrm{M\mathrm{String}}}}}
\def\colim{{{\mathrm{colim}}}}
\def\Spec{{{\mathrm{Spec}}}}
\def\Spf{{{\mathrm{Spf}}}}

\def\Proj{{{\mathrm{Prof}}}}

\def\Aut{{{\mathrm{Aut}}}}
\def\map{{{\mathrm{map}}}}
\def\Gal{{{\mathrm{Gal}}}}
\def\Ga{{{\Gamma}}}

\def\Gr{{{\mathrm{Gr}}}}
\def\RR{{{\mathbb{R}}}}
\def\MO{{{\mathrm{MO}}}}
\def\MU{{{\mathrm{MU}}}}
\def\cH{{{\mathcal{H}}}}
\def\cL{{{\mathcal{L}}}}
\def\cN{{{\mathcal{N}}}}
\def\cU{{{\mathcal{U}}}}
\def\fgl{{{\mathbf{fgl}}}}
\def\Hom{{{\mathrm{Hom}}}}
\def\holim{{\mathrm{holim}}}
\def\defeq{\overset{\mathrm{def}}=}
\def\mm{{{\mathfrak{m}}}}
\def\Def{{{\mathrm{Def}}}}
\def\Proj{{\mathrm{Proj}}}
\def\hatgm{{{\hat{\GG}_m}}}
\def\otherh{C^{\ast}_{\ast}}
\def\grm{{M^\ast_\ast}}

\newcommand{\brackets}[4]{{{\left\{ \begin{array}{ll}
                                {{#1}}&{{#2}}\\
                                \\
                                {{#3}}&{{#4}}\end{array}\right.}}}

\date{Mars 2009}
\bbkannee{61\`eme ann\'ee, 2008-2009}
\bbknumero{1005}
\title{Topological Modular Forms}
\subtitle{after Hopkins, Miller, and Lurie}
\author{Paul G. GOERSS}
\address{Northwestern University\\
Department of Mathematics\\
2033 Sheridan Road\\
EVANSTON, Il 60208-2730 -- USA}
\email{pgoerss@math.northwestern.edu}

\begin{document}
\maketitle

\noindent{\bf INTRODUCTION}

\

In the early 1970s, Quillen \cite{Qui} noticed a strong connection between 
1-parameter formal Lie groups and cohomology theories with a natural theory of
Chern classes. The algebraic geometry of these formal Lie groups allowed Morava, Ravenel,
and others  to make predictions about  large scale phenomena in stable
homotopy theory, and the resulting theorems completely changed the field.
For example, the solution of Ravenel's nilpotence conjectures by Devinatz, 
Hopkins, and Smith (\cite{DHS} and \cite{HS}) was one of the great advances of the 1980s. 

An example of a 1-parameter formal Lie groups can be obtained by taking the formal neighborhood of
the identity in a smooth algebraic group of dimension one. The additive group and the multiplicative
group correspond to ordinary cohomology and complex $K$-theory respectively,
and the only other  algebraic groups of dimension 1 are elliptic curves. 
This class is different because elliptic curves can come in families over 
a base scheme $S$ and the geometry of the fibers can vary significantly as we
move through $S$.  Thus there are many elliptic cohomology theories and it should
be possible to produce them in families over schemes. In retrospect, the realization,
by Hopkins and Miller,  of a good theory of elliptic cohomology theories provided a 
centerpiece for the emerging field of derived algebraic geometry.

Derived algebraic geometry has origins in a number of diverse sources. 
In geometry, there is the work of Serre on multiplicities in intersection theory \cite{Serre} and the work
of Illusie \cite{Ill} on the cotangent complex. For an overview of the roots in
stable homotopy theory, the article \cite{Hop} is very useful. In algebraic
$K$-theory, again originating with Quillen, it was important very early to 
consider algebraic varieties and schemes with sheaves of generalized ring objects and,
indeed, it was mathematicians in this field who first wrote down a systematic
theory \cite{Jar}. Thanks to work of  To\"en, Vezzosi \cite{ToeVez}, and Lurie \cite{Lur},
we now have a fairly mature theory. The purpose here
is to concentrate on the example of elliptic cohomology theories; that is, I would like to make
precise the following statement of a theorem of Mike Hopkins and Haynes Miller,
as refined by Jacob Lurie: the compactified
Deligne-Mumford moduli stack of elliptic curves is canonically and essentially uniquely
an object in derived algebraic geometry. The homotopy global sections of this
derived stack form the ring spectrum of topological modular forms.

\section{AN OVERVIEW}

\subsection{The moduli stack of elliptic curves}

In the late 1960s, Deligne and Mumford \cite{DelMum} defined
a moduli object $\cM_g$ for algebraic curves of genus $g$.  Thus, morphisms $X \to \cM_g$
from a scheme $X$ to $\cM_g$ are in one-to-one correspondence
with smooth proper morphisms 
$$
q: C \longr X
$$
of relative dimension $1$ such that each fiber is a curve of genus $g$. It was known that
that $\cM_g$ could not be a scheme; one way to see this is to note that 
automorphisms of the fibers of $q$ do not vary nicely with the fiber. 
However, Deligne and Mumford noticed that $\cM_g$ exists if we
enlarge the category of schemes slightly to include what we now 
call Deligne-Mumford stacks. From this example, the whole theory of algebraic
stacks emerged. There is an extremely brief exposition on algebraic stacks at the
beginning of section \ref{sectfg}.

From this collection, we single out the moduli stack $\cM_{1,1}$ of
elliptic curves. These are curves of genus $1$ with $1$ marked point;
that is diagrams of the form
\[
\xymatrix{
C \ar@<.5ex>[r]^q &  \ar@<.5ex>[l]^e X
}
\]
where $q$ is a family of smooth curves of genus $1$ over $X$ and 
$e$ is a section identifying a distinguished point in each fiber. A classical,
but still remarkable, property of these curves is that  $C$ becomes
an abelian group over $X$: there is a canonical commutative group multiplication
on $C$
with $e$ as the unit. Furthermore, any morphism of elliptic curves is a group
homomorphism. There are a number of reasons for singling out this
stack: elliptic curves are central to algebraic number theory, for example. 
For algebraic topologists, the formal neighborhood $C_e$ of $e$ in
$C$ gives a family of 1-parameter formal Lie groups (``formal groups'')
which, in turn, gives rise to a rich family of elliptic cohomology theories.

The stack $\cM_{1,1}$ is not compact -- more precisely, it is not proper
over $\ZZ$. Indeed, the morphism to the affine line
$$
j:\cM_{1,1} \longr \AA^1
$$
which assigns to each elliptic curve $C$ the $j$-invariant $j(C)$
is almost (but not quite) a covering map. Deligne and Mumford found
a canonical compactification $\mell$ of $\cM_{1,1}$
which now classifies generalized elliptic curves with, possibly, nodal singularities.
In an extensive study of this stack \cite{DelRap}
Deligne and Rappoport showed that the $j$-invariant extends to
a morphism
$$
j:\mell \longr \PP^1.
$$

\subsection{Derived Schemes}

The basic objects of algebraic geometry are schemes, which are locally
ringed spaces $(X,\cO_X)$ which are locally isomorphic to the prime
spectrum of some ring. Thus, among other things, $X$ is a space and
$\cO_X$ is a sheaf of rings on $X$. The basic idea of derived algebraic
geometry is to replace rings by some more generalized ring object. 
For example, Serre \cite{Serre} considered schemes
with a sheaf of commutative differential graded algebras. This has the
advantage that it's relatively simple to define and, indeed, over the rational numbers
it is equivalent to  the more general theory. However, commutative DGAs only have
good homotopy theory when we work over the rational numbers; over the
integers or in characteristic $p$ a more rigid theory is needed. In his work
on the cotangent complex, Illusie \cite{Ill} worked with schemes with sheaves
of simplicial commutative algebras. This was also the point of view of Lurie
in his thesis \cite{Lurthesis}. However, there are basic examples arising from homotopy theory
which cannot come from simplicial algebras -- complex $K$-theory is
an important example. Thus, a derived scheme (or stack) will be a scheme
equipped with a sheaf of commutative ring spectra. I immediately remark
that ``commutative'' is a difficult notion to define in ring spectra:
here I mean `$`E_\infty$-ring spectra''. 
The foundations of commutative ring spectra
are forbidding, but I'll make some attempt at an exposition below
in Section \ref{spectra}.

A spectrum $X$ has homotopy groups $\pi_kX$, for $k \in \ZZ$. If $X$ is a
commutative ring spectrum, $\pi_0X$ is a commutative ring and the graded
abelian group $\pi_\ast X$ is a graded skew-commutative $\pi_0X$-algbera. In
particular, $\pi_kX$ is a $\pi_0X$-module. 

\begin{defi}
\label{derived-schm} A {\bf derived scheme} $(X,\cO)$
is a pair with $X$ a topological space and $\cO$ a sheaf of commutative
ring spectra on $X$ so that
\begin{enumerate}

\item the pair $(X,\pi_0\cO)$ is a scheme; and

\item the sheaf $\pi_k\cO$ is a quasi-coherent sheaf of $\pi_0\cO$-modules.
\end{enumerate}
\end{defi}
It is somewhat subtle to define the notion of a sheaf of spectra; I will come back
to that point below in Remark \ref{sheaf-sp}.  One definition of ``quasi-coherent'' 
is to require the sheaf to be locally the cokernel of a morphism
between locally free sheaves. On the spectrum of a ring $R$, these are the module
sheaves which arise from $R$-modules. There is a completely analogous definition of a derived
Deligne-Mumford stack, except that now we must be careful about the topology
we use to define sheaves: for these objects we must use the \'etale topology.

There is rich structure inherent in Definition \ref{derived-schm}. The homotopy
groups of a commutative ring spectrum support far more structure than simply that of
a graded skew-commutative ring; in particular, it is a ring with ``power operations''.
See Remark \ref{after-LMS} below. Thus, if $(X,\cO)$
is a derived scheme, the graded sheaf $\pi_\ast \cO$ is a sheaf of graded rings with
all of this higher order structure.

\subsection{Topological modular forms}

On the compactified Deligne-Mumford stack $\mell$, there is a canonical
quasi-coherent sheaf $\omega$. If $C$ is an elliptic curve over $X$,
then $C$ is an abelian variety over $X$ of relative dimension $1$ and we can construct the 
sheaf of invariant $1$-forms $\omega_C$ for $C$. This is locally free sheaf
of rank 1 on $X$ and the assignment
$$
\omega(C:X \to \mell) = \omega_C
$$
defines a quasi-coherent sheaf on $\mell$. The sheaf $\omega$ is a locally
free of rank 1, hence invertible, and the tensor powers $\omega^{\otimes k}$,
$k \in \ZZ$ are all quasi-coherent. Here is the main result; see
\cite{Hop1} and \cite{Lurell}.

\begin{theo}[{\bf Hopkins-Miller-Lurie}]
\label{mn-thm}
There exists a derived Deligne-Mumford stack $(\mell,\cO)$ so that
\begin{enumerate}

\item the underlying algebraic stack
$(\mell,\pi_0\cO)$ is equivalent to the compactified Deligne-Mumford
moduli stack $(\mell,\cO_\Ell)$ of generalized elliptic curves; and

\item there are isomorphisms of quasi-coherent sheave $\pi_{2k}\cO \cong
\omega^{\otimes k}$ and $\pi_{2k+1}\cO = 0$.
\end{enumerate}
Furthemore, the derived stack $(\mell,\cO)$ is determined up to equivalence by 
conditions (1) and (2).
\end{theo}

The doubling of degrees in (2) is quite typical: when the homotopy groups
of a commutative ring spectrum $X$ are concentrated in even degrees, then
$\pi_\ast X$ is a graded commutative ring, not just commutative up to sign.

We can rephrase the uniqueness statement of Theorem \ref{mn-thm} as follows: there
is a space of all derived stacks which satisfy points (1) and (2) and this
space is path-connected. In fact, Lurie's construction gives a canonical
base-point; that is, a canonical model for $(\mell,\cO)$.

Here is a definition of the object in the title of this manuscript. 

\begin{defi}
\label{tmf-defn}
The spectrum $\tmf$ of {\bf topological modular forms} is the derived
global sections of the sheaf $\cO$ of commutative ring spectra
of the derived Deligne-Mumford stack $(\mell,\cO)$. 
\end{defi}

The spectrum $\tmf$ is itself a commutative ring spectrum and 
formal considerations now give a descent spectral sequence
\begin{equation}
\label{main-desc}
H^{s}(\mell,\omega^{\otimes t}) \Longrightarrow \pi_{2t-s}\tmf.
\end{equation}
This spectral sequence has been completely calculated in \cite{Bau}, 
\cite{HopMah}, and \cite{Rez}.

The spectrum $\tmf$ is called ``topological modular forms'' for the following
reason. One definition of modular forms of weight $k$ (and level $1$) is
as the global sections of $\omega^{\otimes k}$ over $\mell$;
that is, the graded ring
$$
M_\ast = H^0(\mell,\omega^{\otimes \ast})
$$
is the ring of modular forms for generalized elliptic curves. This
ring is well understood \cite{Del}: there are modular forms
$c_4$, $c_6$, and $\Delta$ of weights $4$, $6$, and $12$ respectively
and an isomorphism of graded rings
$$
\ZZ[c_4,c_6,\Delta]/(c_4^3 - c_6^2 = 1728\Delta) \cong M_\ast.
$$
Note $1728 = (12)^3$, indicating that the primes $2$ and $3$ are special
in this subject. The modular form $\Delta$ is the discriminant and is the test for
smoothness: a generalized elliptic curve $C$ is smooth if and only
if $\Delta(C)$ is invertible; in fact, $\cM_{1,1} \subseteq \mell$ is the
open substack obtained by inverting $\Delta$. 

Modular forms then form the zero line of the spectral sequence
of (\ref{main-desc}), at least up to degree doubling,
and we can now ask which modular forms
give homotopy classes in $\tmf$. The higher cohomology groups
of (\ref{main-desc}) are all $2$ and $3$-torsion and, as mentioned above,
the differentials have also been calculated. Thus, for example,
we know that $\Delta$ is not a homotopy class, but $24\Delta$ is;
similarly $c_6$ is not a homotopy class, but $2c_6$ is. The class
$c_4$ is a homotopy class. In the last section of this note, I will
uncover some of the details of this calculation.

\subsection{Impact}

Let me write down three areas of algebraic topology where $\tmf$
has had significant impact. I'll go in more-or-less chronological order.

\begin{rema}[{\bf The Witten genus}]\label{witten} In his work in
string theory, Witten noticed that one could define a
genus for certain spin manifolds that takes values in modular
forms. A compact differentiable manifold $M$ has a spin structure if 
the first and second Stiefel-Whitney classes vanish and
a spin manifold is $M$ is a manifold with a choice of spin
structure. For spin manifolds there is a new characteristic class
$\lambda$ which has the property that
twice this class is the first Pontrjagin class. If $\lambda$ also
vanishes, the manifold $M$ has a {\bf string} structure. 
There is a cobordism ring $\mstring_\ast$ of string manifolds,
and Witten wrote down a rich ring homomorphism
\[
\sigma_\ast:\mstring_\ast \longr \QQ\otimes M_\ast,
\]
which we now call the Witten genus. Let me say a little about his
methods.

Every modular form can be written as a power series over the integers
in $q$; this is the $q$-expansion. One way to do this is to evaluate the modular 
form on the Tate curve, which is a generalized elliptic curve over $\ZZ[[q]]$. 
This defines a monomorphism $M_\ast \to \ZZ[[q]]$ but it is definitely not
onto. For each string manifold, Witten \cite{Wit} wrote down a power series
in $q$ over the rationals and then used physics to argue it must be a modular
form. 

Almost immediately, homotopy theorists began searching for a spectrum-level
construction of this map. By the Thom-Pontrjagin construction,
$\mstring_\ast$ is isomorphic to the homotopy groups of a commutative
ring spectrum $\mstring$ and the Witten genus was posited to be given by a morphism
of ring spectra from $\mstring$ to some appropriate ring spectrum or family of
ring spectra; this is  what happens with 
the Atiyah-Bott-Shapiro realization of the $\hat{A}$-genus on spin
manifolds. Motivated by this problem, among others, there was important early important
work on elliptic cohomology theories and orientations and a major conference on 
elliptic cohomology theories in the late
1980s;  \cite{LRS} is a highlight of
this period. There were also very influential papers in the years
following. For example, the Bourbaki expos\'e by  Segal \cite{Seg} has been
especially important in the search for elliptic cohomologies
arising from differential geometry and mathematical physics and the paper
by Franke \cite{Fra} gave a profusion of examples of cohomology
theories arising from elliptic curves.

The Witten genus was addressed specifically
in the work of Ando, Hopkins, Rezk, and Strickland. Earlier results are in
 \cite{AHS} and \cite{Hop2}; the definitive results are in
\cite{AHR1} and \cite{AHR2}.

\begin{theo}\label{ahr}The Witten genus
can be realized as a morphism of commutative ring spectra
\[
\sigma:\mstring \longr \tmf.
\]
The map $\sigma$ is surjective on homotopy in non-negative degrees.
\end{theo}

Besides giving a rigid construction of the Witten genus, this has other
consequences. Not every modular form is in the image, for example,
as not every modular form is a homotopy class. On the other hand,
there is a plenty of $2$ and $3$ torsion in both source and target
and the map $\sigma$ detects whole families of this torsion.
\end{rema}

\begin{rema}[{\bf Homotopy groups of spheres}]\label{homotopy} 
One of the fundamental problems of stable homotopy theory is to compute $\pi_\ast S$,
where $S$ is the stable sphere; thus $\pi_k S = \colim \pi_{n+k}S^n$.
This is impossible, at least at this point, but
we can measure progress against this problem. Early successful calculation focused
on the image of the $J$-homomorphism $\pi_\ast SO(n) \to \pi_\ast S^n$; see \cite{Ada1} 
and \cite{Mah}. There were many calculations with higher order phenomena,
most notably in \cite{MRW} and the work of Shimomura and his coauthors (see \cite{Shi} and \cite{ShiWan} among
many papers), but  the very richness of these results hindered comprehension.
Using $\tmf$, related spectra, and the algebra and geometry of ellipitic curves, it is
now possible to reorganize the calculations in a way that better reveals the larger
structure. See for example, \cite{GHMR}, \cite{HKM}, and \cite{Beh1}. The latter,
in particular, makes the connections with elliptic curves explicit.
\end{rema}

\begin{rema}[{\bf Congruences among modular forms}]\label{congruences}
There is a remarkable interplay between homotopy theory and the theory
of modular forms. For example, by thinking about the descent spectral sequence
of Equation (\ref{main-desc}) above, Hopkins \cite{Hop1} discovered a new proof
of a congruence of Borcherds \cite{Bor}. We can also vary the moduli problem
to consider elliptic curves with appropriate homomorphisms $(\ZZ/N\ZZ)^2 \to C$
(``level-structures''). This new moduli problem is \'etale over $\cM_{1,1}[1/N]$
and Theorem \ref{mn-thm} immediately produces a new spectrum $\tmf[N]$ in
$\ZZ[1/N]$-local homotopy theory, which might be called topological modular forms
of level $N$. There is a descent spectral sequence analogous to  (\ref{main-desc})
beginning with modular forms
of level $N$. This has been used to effect in \cite{MR2} and \cite{Beh1}. 

More subtle connections and congruences emerge when we work $p$-adically.
For example,  the ring of divided congruences of Katz \cite{Kat} appears
in topology as the $p$-complete $K$-theory of $\tmf$. This observation, due to 
Hopkins, appears in \cite{Lau} for the prime $2$. More recently,  Behrens \cite{Beh2} was
able to explicitly compare certain congruences among modular forms that hold
near supersingular curves with the intricate and beautiful patterns in the stable homotopy groups
of spheres that first appeared in \cite{MRW}.
\end{rema}

\section{Basics from homotopy theory}

\subsection{Spectra and commutative ring spectra}\label{spectra}

The need for a good theory of spectra arose in the 1950s while trying to find
a framework to encode a variety of examples that displayed similar phenomena.
In each case there was a {\bf spectrum}; that is, a sequence of based (or pointed) spaces $X_n$ and suspension
maps
$$
\Sigma X = S^1 \wedge X_n \to X_{n+1}.
$$
Here $A \wedge B$ is the smash product of two pointed spaces characterized
by the adjunction formula in spaces of pointed maps
$$
\map_\ast (A \wedge B,Y) \cong \map_\ast (A,\map_\ast(B,Y)).
$$ 

\begin{exem}\label{the-basic-exam}
1.) The space $S^1 \wedge A$ is the suspension of $A$ and a basic example
is the suspension spectrum $\Sigma^\infty A$ of a pointed finite
CW complex $A$; thus, $(\Sigma^\infty A)_n = \Sigma^n A$ is the iterated suspension.
In particular, $S = \Sigma^\infty S^0$.

2.) A generalized cohomology theory is a contravariant functor $E^\ast(-)$ from the
homotopy category of spaces to graded abelian groups which takes disjoint unions
to products and has a Mayer-Vietoris sequence. By Brown's representability theorem
there are spaces $E_n$ so the reduced cohomology $\tilde{E}^n(A)$ is naturally
isomorphic to  the pointed homotopy classes of maps $A \to E_n$.
The suspension
isomorphism $\tilde{E}^{n}A \cong \tilde{E}^{n+1}\Sigma A$ yields the map
$\Sigma E_n \to E_{n+1}$.

3.) The Thom spectra $\MO$ and its variants such as $\mstring$. The $n$th space of $\MO$ is
the Thom space $\MO(n)=T(\gamma_n)$
of the universal real $n$-plane bundle $\gamma_n$ over the infinite Grassmannian.
The map $\Gr_n(\RR^\infty) \to \Gr_{n+1}(\RR^\infty)$ classifying the Whitney sum
$\epsilon\oplus \gamma_n$, where $\epsilon$ is a trivial $1$-plane bundle,
defines the map
$$
S^1 \wedge \MO(n) \cong T(\epsilon \oplus \gamma_n) \to T(\gamma_{n+1}) = \MO(n+1).
$$
\end{exem}

A basic invariant of a spectrum $X$ is the {\bf stable homotopy groups}
$$
\pi_k X = \colim\ \pi_{n+k}X_n,\qquad k \in \ZZ.
$$
For example, the Thom-Pontrjagin construction shows that $\pi_k \MO$ is the 
group of cobordism classes of closed differentiable $k$-manifolds; it was
an early triumph of homotopy theory that Thom was able to compute these groups.

\begin{rema}[{\bf The stable homotopy category}]\label{stable-homotopy}
A basic difficulty was to decide how to build a homotopy theory of spectra.
Seemingly anomalous examples led to some very ingenious ideas,
most notably the ``cells now maps later'' construction of Boardman and 
Adams \cite{AdamsBlue}. The language
of model categories gives a simple and elegant definition. We define a morphism
$f:X \to Y$ of spectra to be a set of pointed maps $f_n:X_n \to Y_n$ which commute with the suspension
maps. Next
we stipulate that such a morphism is a weak equivalence if
$\pi_\ast f: \pi_\ast X \cong \pi_\ast Y$. Bousfield and Friedlander \cite{BF} showed that there
is a model category structure on spectra with these weak equivalences, and we obtain
the homotopy category by formally inverting the weak equivalences.
\end{rema}

There was a much more serious difficulty in the theory, however.  Suppose that we have a cohomology theory
$E^\ast(-)$ with natural, associative, and graded commutative cup products
$$
E^mX \otimes E^n X \longr E^{m+n}X.
$$
Then, by considering the universal examples, we get a map
$$
E_m \wedge E_n \longr E_{n+m}.
$$
From this we ought to get a map $E \wedge E \to E$ of spectra making $E$ into
a homotopy associative and commutative ring object. However, given two
spectra $X$ and $Y$ the object $\{X_m \wedge Y_n\}$ isn't a spectrum, but
a ``bispectrum''. We can extract a spectrum $Z$ by taking $Z_k$ to be any
of the spaces $X_m \wedge Y_n$ with $m+n=k$ with $m$ and $n$ non-decreasing; if, in addition, we
ask that
$n$ and $m$ go to $\infty$ with $k$, we get a well-defined homotopy type --
and hence a symmetric monoidal smash product on the homotopy category. 
However,  this does not descend from such a structure on the category of
spectra.

Worse, there were important examples that indicated that there should be such
a structure on spectra. The Eilenberg-MacLane spectra $HR$, $R$ a commutative ring, the
cobordism spectrum $\MO$ and the complex analog $\MU$, spectra arising
from algebraic and topological $K$-theory, and the sphere spectrum $S$ itself
all had more structure than simply giving ring objects in the homotopy category.

\begin{rema}\label{LMS}The first, and still a very elegant, solution to this problem was due to Lewis,
May, and Steinberger \cite{LMS}. The idea was to expand the notion
of a spectrum $X$ to be a collection of spaces $\{X_V\}$ indexed on the finite-dimensional
subspaces of the infinite inner product space $\RR^\infty$; the suspension maps
then went $S^V \wedge X_W \to X_{V\oplus W}$, where $S^V$ was the one-point
compactification of $V$. If $X_i$, $1 \leq i \leq n$ are spectra, then $X_1 \wedge\cdots \wedge X_n$
naturally yields a spectrum object over
$(\RR^\infty)^n$ and any linear isometry $f:\RR^\infty \to (\RR^\infty)^n$
 then returns a spectrum $f^\ast(X_1 \wedge \cdots \wedge X_n)$ over 
 $\RR^\infty$. The crucial observation is that the
space $\cL(n)$ of all such choices of isometries is contractible and the construction
of $f^\ast(X_1 \wedge \cdots \wedge X_n)$
could be extended to a construction of a functor
$$
(X_1,\cdots,X_n) \mapsto \cL(n)_+ \wedge X_1 \wedge \cdots \wedge X_n.
$$
This functor, in the case $n=2$, descended to the smash product on the homotopy category. 
More importantly, the collection of spaces $\cL = \{\cL(n)\}$ form an $E_\infty$-operad. This
implies that there is a coherent way to compose the above functors and that  the action
of the symmetric group on $\cL(n)$ is free. 
\end{rema}

We still didn't have a spectrum level symmetric monoidal structure, but we did now know
what a ring object should be.

\begin{defi}\label{comm-ring-spec} A {\bf commutative ring spectrum} is an algebra over the 
operad $\cL$; that is, an algebra for the monad
$$
X \mapsto \vee_n\ \cL(n)_+ \wedge_{\Sigma_n} X^{\wedge n}.
$$
Here $\vee$ is the coproduct or wedge in spectra; and $\wedge_{\Sigma_n}$ means
divide out by the diagonal action of the symmetric group. We can also call a commutative
ring spectrum an $E_\infty$-ring spectrum.
\end{defi}

This definition has the distinct advantage
of being a machine with an input slot; for example, this theory is ideal for showing that
the Thom spectra $\MO$ and $\MU$ are commutative ring spectra;
similarly, the Eilenberg-MacLane spectra $HR$, with $R$ a commutative ring, and
$S$ are easily seen to be commutative ring spectra.

\begin{exem}\label{absorb} If $R$ is a commutative ring, we could
consider algebras over an $E_\infty$-operad in $R$-chain complexes;
that is, we could work with differential graded $E_\infty$-algebras. However, in a very
strong sense (the technical notion is {\it Quillen equivalence} or equivalence of
$\infty$-categories) the
category of differential graded $E_\infty$-algebras over $R$  has the same homotopy theory
as the category of commutative
$HR$-algebras in spectra, where $HR$ is the associated Eilenberg-MacLane
spectra. This is actually a fairly restrictive example of commutative ring spectra; many important
examples, such as $\MO$, $\MU$, and $K$-theory are not of this type.
I note that  if $R$ is a $\QQ$-algebra, then $E_\infty$-algebras over $R$ and commutative
dgas over $R$ also form Quillen equivalent categories.
\end{exem}

\begin{rema}\label{after-LMS}
Since the work of Lewis, May, and Steinberger, many authors (including May himself) have
built models for the stable category which indeed have a symmetric monoidal structure. See 
\cite{EKMM}, \cite{HSS}, and \cite{MMSS}. Then
a commutative ring spectrum is simply a commutative monoid for that structure. However, it
is only a slight
exaggeration to say that all such models build an $E_\infty$-operad into the spectrum-level
smash product in some way or another; thus, the theory is elegant, but the computations remain
the same.  For example, the homotopy of  an $E_\infty$-ring spectrum is a commutative
ring, but there is much more structure as well: the map in homotopy
$$
\pi_\ast \cL(n)_+ \wedge_{\Sigma_n} X^{\wedge n} \longr \pi_\ast X
$$
adds {\bf power operations} (such as Steenrod or Dyer-Lashof operations) to $\pi_\ast X$.
The spectral sequence
$$
H_p(\Sigma_n, \pi_q( X^{\wedge n})) \Longrightarrow
\pi_{p+q} \cL(n)_+ \wedge_{\Sigma_n} X^{\wedge n}
$$
makes the role of the homology of the symmetric groups explicit in the construction these
operations.
\end{rema}

\begin{rema}[{\bf Sheaves of spectra}]\label{sheaf-sp} Once we have ring spectra, we must
confront what we mean by a sheaf of ring spectra. There is an issue here as well. Suppose
we have a presheaf $\cF$ on $X$ and $\{V_i\}$ is a cover of $U\subseteq X$. Then if $\cF$ is a sheaf we have an
equalizer diagram
$$
\xymatrix{
\cF(U) \ar[r] & \prod \cF(V_i) \ar@<.5ex>[r] \ar@<-.5ex>[r] & \prod \cF(V_i \times_U V_j).
}
$$
The problem is that equalizers are not homotopy invariant. Again model categories help.
We define a morphism of presheaves $\cE \to \cF$ to be a weak equivalence
if it induces an isomorphism $\pi_\ast \cE \to \pi_\ast \cF$ of associated homotopy sheaves. 
Thus, for example, a presheaf is weakly equivalent to its associated sheaf. The theorem,
due to Jardine \cite{Jar}, \cite{JarR}, is that there is a model category structure on presheaves with these weak
equivalences. A {\bf sheaf of (ring) spectra} is then a fibrant/cofibrant object in this model
category. This has the effect of building in the usual homological algebra for sheaves;
that is, if $\cF$ is a module sheaf and $\cF \to \cI^\bullet$ is an injective resolution; then
the associated presheaf of generalized Eilenberg-MacLane spectra $K\cI^\bullet$ is a
fibrant model for $K\cF$. In this setting, global sections $\Ga(-)
$ are inherently derived so I may write $R\Ga(-)$ instead. In good cases there is a descent spectral sequence
$$
H^s(X,\pi_t\cF) \Longrightarrow \pi_{t-s}\Ga(\cF).
$$ 
\end{rema}

\subsection{Cohomology theories and formal groups}

If $E$ is a spectrum associated to a cohomology theory $E^\ast$, we get a homology
theory by setting $E_\ast X = \pi_\ast E \wedge X_+$, where $X_+$ is $X$ with
a disjoint basepoint.

\begin{defi}\label{2-period} Let $E^\ast(-)$ be a cohomology theory. Then $E^\ast$
is {\bf 2-periodic} if
\begin{enumerate}

\item the functor $X \mapsto E^\ast(X)$ is a functor to graded commutative rings;

\item for all integers $k$, $E^{2k+1}=E^{2k+1}(\ast) = 0$;

\item $E^2$ is a projective module of rank $1$ over $E^0$; and

\item for all integers $k$, the cup product map $(E^2)^{\otimes k} \to E^{2k}$
is an isomorphism.
\end{enumerate}
\end{defi}

Note that  $E^2$ is an invertible module over $E^0$ and $E^{-2}$ is the dual
module.
If $E^2$ is actually free, then so is $E_2 = E^{-2}$ and a choice of generator
$u \in E_2$ defines an isomorphism $E_0[u^{\pm 1}] \cong E_\ast$. 
This often happens; for example,
in complex $K$-theory. However,
there are elliptic cohomology theories for which $E_2$ does not have a global
generator.

From $2$-periodic cohomology theories we automatically get a formal group. Let
$\CP^\infty = \Gr_1(\CC^\infty)$
be the infinite Grassmannian classifying complex line bundles. Then 
$\CP^\infty$ is a topological monoid where the multiplication $\CP^\infty \times
\CP^\infty\to \CP^\infty$ classifies the tensor product of line bundles.
If $E^\ast$ is a $2$-periodic homology theory, then $E^0\CP^\infty$ is
complete with respect to the augmentation ideal 
$$
I(e) \defeq \tilde{E}^0 \CP^\infty = \mathrm{Ker}\{\ E^0 \CP^\infty \to E^0(\ast)\ \}
$$
and,  using the monoid structure on $\CP^\infty$, we get a commutative
group object in formal schemes
$$
G_E = \Spf(E^0 \CP^\infty).
$$
This formal group is smooth and one-dimensional in the following sense. Define
the $E_0$-module $\omega_G$ by
\begin{equation}\label{inv-diff-defn}
\omega_G = I(e)/I(e)^2 \cong \tilde{E}^0S^2 \cong E_2.
\end{equation}
This module is locally free of rank 1, hence projective, and any choice of
splitting of $I(e) \to \omega_G$ defines a homomorphism out of the
symmetric algebra
$$
S_{E_0}(\omega_G) \longr E^0\CP^\infty 
$$
which becomes an isomorphism after completion.
For example, if $E_2$ is actually free
we get a non-canonical isomorphism
$$
E^0\CP^\infty \cong E^0[[x]].
$$
Such an $x$ is called a {\it coordinate}. Whether $G_E$ has a coordinate or not, the
ring $E_0$ and the formal group determine the graded coefficient ring $E_\ast$; indeed,
$E_{2t+1}=0$ and for all $t \in \ZZ$,
\begin{equation}\label{inv-diff-defn-ext}
\omega_G^{\otimes t}  \cong \tilde{E}^0S^{2t} \cong E_{2t}.
\end{equation}

\begin{rema}[{\bf Formal group laws}]\label{why-fgs-not-fgls}The standard literature on chromatic
homotopy theory, such as \cite{AdamsBlue} and \cite{Rav}, emphasizes formal
group laws. If $E^\ast(-)$ is a two-periodic theory with a coordinate,
then the group multiplication
$$
G_E \times G_E = \Spf(E^0(\CP^\infty \times \CP^\infty)) \to
\Spf(E^0 \CP^\infty) = G_E
$$
determines and is determined by a power series
$$
x+_Fy = F(x,y) \in E^0[[x,y]] \cong E^0(\CP^\infty \times \CP^\infty).
$$
This power series is a $1$-dimensional formal group law. 

Homomorphisms can also be described by power series, and a homomorphism
$\phi$ is an isomorphism if $\phi'(0)$ is a unit.
\end{rema}

\begin{rema}[{\bf Invariant differentials}]\label{invariant differential} The module
$\omega_{G_E}$ is defined as the conormal module of the embedding
$$
e:\Spec(E^0) \to \Spf(E^0\CP^\infty) = G_E
$$
given by the basepoint.  This definition extends to any formal group $G$ over
any base scheme $X$ and we get a sheaf $\omega_G$ on $X$. It is
isomorphic to the sheaf of invariant differentials of $G$ via an evident
inclusion $\omega_G \to q_\ast \Omega_{G/X}$. Here
$\Omega_{G/X}$ is the sheaf of continuous differentials and $q:G \to X$
is the structure map.

If $X = \Spec(R)$ for some ring and $G$ has a coordinate $x$, then the invariant differentials
form the free $R$-module generated
by the {\it canonical} invariant differential
$$
\eta_G = \frac{dx}{F_y(x,0)}
$$
where $F_y(x,y)$ is the partial derivative of the associated formal group
law. It is an exercise to calculate that if $\phi:G_1 \to G_2$ is a homomorphism
of formal groups with coordinate, then
$d\phi:{\omega}_{G_2} \to {\omega}_{G_1}$
is determined by
\begin{equation}\label{inv2}
d\phi(\eta_{G_2}) = \phi'(0)\eta_{G_1}.
\end{equation}
\end{rema}

\section{Formal groups and stable homotopy theory}

\subsection{The moduli stack of formal groups}\label{sectfg}

Let $\cM_\fg$ be the moduli stack of formal groups: this is the algebro-geometric
object which classifies all smooth 1-parameter formal Lie groups and their
isomorphisms. Thus, if $R$ is a commutative ring, the morphisms
$$
G:\Spec(R) \longr \cM_\fg
$$
are in one-to-one correspondence with formal groups $G$ over $R$.
Furthermore, the 2-commutative diagrams
\begin{equation}\label{2-commute}
\xymatrix@R=10pt{
\Spec(S) \ar[dr]^H \ar[dd]_f\\
&\cM_\fg\\
\Spec(R)\ar[ur]_G 
}
\end{equation}
correspond to pairs $(f:R\to S,\xymatrix{\phi:H \ar[r]^{\cong}& f^\ast G})$.

\begin{rema}\label{stack-scheme-functor} Schemes are 
defined as locally ringed spaces $(X,\cO_X)$ which have
an open cover, as locally ringed spaces, by affine schemes. 
Equivalently, schemes can be defined as functors from rings to 
sets which are sheaves in the Zariski topology and have an
open cover, as functors, by functors of the form
$$
A \mapsto \mathbf{Rings}(R,A).
$$
Stacks are generalization of the the second definition. A stack is 
a sheaf of groupoids on commutative rings satisfying an  
effective descent condition \cite{Laumon}\S 3. For example, $\cM_\fg$ assigns to 
each ring $R$ the groupoid of  formal
groups over $\Spec(R)$.\footnote{As in \cite{Laumon},
\S 2,  we should really
speak of categories fibered in groupoids, rather than sheaves of
groupoids -- for $f^\ast g^\ast G$ is only isomorphic to $(gf)^\ast G$. However,
there are standard ways to pass between the two notions.}

Algebraic stacks have a suitable cover by schemes.
A morphism $\cM \to \cN$ of stacks is {\it representable}
if for all morphisms $X \to \cN$ with $X$ a scheme, the $2$-category
pull-back (or homotopy pull-back) $X \times_\cN \cM$
is equivalent to a scheme. A representable morphism then has algebraic
property $P$ (flat, smooth, surjective, \'etale, etc.) if all the resulting morphisms
$$
X \times_\cN \cM \to X
$$
have that property.

A stack $\cM$ is then called {\it algebraic}\footnote{The notion defined here is stronger than
what is usually called an algebraic (or Artin) stack, which requires a cover only by an algebraic
space. Algebraic spaces are sheaves which themselves have an appropriate cover by a scheme. Details are in
\cite{Laumon}.} if
\begin{enumerate}

\item every morphism $Y \to \cM$ with $Y$ a scheme is representable; and

\item there is a smooth surjective map $q:X \to \cM$ with $X$ a scheme. 
\end{enumerate}

The morphism $q$ is called a presentation. Note that an algebraic
stack may have many presentations. If a presentation can be
chosen to be \'etale, we have a {\it Deligne-Mumford stack}.
\end{rema}

\begin{rema}\label{not-alg-stack}
The stack $\cM_\fg$ is not algebraic, in this sense, as it only has
a flat presentation, not a smooth presentation. If we define $\fgl$ to
be the functor which assigns to each ring $R$ the set of formal group laws
over $R$, then Lazard's theorem \cite{Laz} says that $\fgl = \Spec(L)$
where $L$ is (non-canonically) isomorphic to $\ZZ[t_1,t_2,\ldots]$.
The map
$$
\fgl \longr \cM_\fg
$$
which assigns a formal group law to its underlying formal group is flat
and surjective, but not smooth since it's not finitely presented.
It is {\it pro-algebraic} however; that is, it can be 
written as the 2-category inverse limit of a tower of the algebraic
stacks of ``buds'' of formal groups.
This is inherent in \cite{Laz} and explicit in \cite{Smith}.
\end{rema}

\begin{rema}\label{mod-sheaves} A sheaf in the $fpqc$-topology
on an algebraic stack $\cM$ is a functor $\cF$ on the category of affine
schemes over $\cM$ which satisfies faithfully flat
descent.  For example, define the structure sheaf $\cO_\fg$
to be the functor on affine schemes over $\cM_\fg$ with
$$
\cO_\fg(R,G) = \cO_\fg(G:\Spec(R) \to \cM) = R.
$$
A module sheaf $\cF$ over $\cO_\cM$ is {\it quasi-coherent} if, for each $2$-commutative
diagram over $\cM$, the restriction map $\cF(R,G) \to \cF(S,H)$ extends to
an isomorphism
$$
S \otimes_R \cF(R,G) \cong \cF(S,H).
$$
This isomorphism can be very non-trivial, as it depends on the 
choice of isomorphism $\phi$ which makes the diagram
$2$-commute.

A fundamental example of a quasi-coherent sheaf is the sheaf 
of invariant differentials $\omega$ on $\cM_\fg$
with
$$
\omega(R,G) = \omega_G
$$
the invariant differentials on $G$. This is locally free of rank $1$ and hence
all powers $\omega^{\otimes n}$, $n \in \ZZ$, are also quasi-coherent sheaves.
The effect of the choice of isomorphism in the $2$-commuting diagram on the
transition maps for $\omega^{\otimes n}$ is displayed in Equation
\ref{inv2}.
\end{rema}

\subsection{The height filtration}

Consider a formal group $G$ over a scheme $X$ over $\FF_p$.
If we let $f:X \to X$ be the Frobenius, we get a new formal group
$G^{(p)} = f^\ast G$. The Frobenius $f:G \to G$
factors through the {\it relative Frobenius} $F:G \to G^{(p)}$.
We know that if $\phi:G \to H$
is a homomorphism of formal groups over $X$
for which $d\phi = 0:\omega_H \to \omega_G$, there is then a factoring
$\phi = \psi F: G \to H$.
Then we can test $d\psi$ to see if we can factor further.

For example, let $\phi = p:G \to G$ be $p$th power map. Then we obtain
a factoring 
$$
\xymatrix{
G \ar[r]_-F  \ar@/^1pc/[rr]^p& G^{(p)} \ar[r]_-{V_1} &G
}
$$
This yields an element
$$
dV_1 \in \Hom(\omega_G,\omega_{G^{(p)}})
$$
and we can factor further if $dV_1=0$. Since $G$ is of dimension 1,
$\omega_{G^{(p)}}= \omega_G^{\otimes p}$; since $\omega_G$ is invertible,
$$
 \Hom(\omega_G,\omega_G^{\otimes p}) = \Hom(\cO_X,\omega_G^{\otimes p-1}).
$$
Thus $dV_1$ defines a global section
$v_1(G)$ of $\omega_G^{\otimes p-1}$. If $v_1(G) = 0$, then we obtain
a further factorization and a global section 
$v_2(G) \in \omega_{G}^{\otimes p^2-1}$. This can be continued
to define sections $v_n(G) \in \omega_G^{\otimes p^n-1}$ and $G$ has {\bf height}
at least $n$ if
$$
v_1(G) = \cdots = v_{n-1}(G) = 0.
$$
We say $G$ has height
{\it exactly} $n$ if $v_n(G):\cO_X \to \omega_G^{\otimes p^n-1}$ is an isomorphism.
Note that a formal group may have infinite height.

The assignment $G \mapsto v_1(G)$
defines a global section $v_1$ of the sheaf $\omega^{\otimes p-1}$ on
the closed substack
$$
\FF_p \otimes \cM_\fg \defeq \cM(1) \subseteq \cM_\fg
$$
Indeed, we obtain a sequence of closed substacks
$$
\dots \subseteq \cM(n+1) \subseteq \cM(n) \subseteq
\cdots \subseteq \cM(1) \subseteq \cM_\fg
$$
where $\cM(n+1) \subseteq \cM(n)$ is defined by the vanishing of the
global section $v_n$ of $\omega^{\otimes p^n-1}$. Thus $\cM(n)$ classifies
formal groups of height at least $n$. The relative open
$$
\cH(n) = \cM(n) - \cM(n+1)
$$
classifies formal groups of height exactly $n$. One of Lazard's theorems \cite{Laz}, rephrased,
says that $\cH(n)$ has a single geometric point given by a formal group
$G$ of height $n$ over any algebraically closed field $\FF$ of characteristic $p$. The pair
$(\FF,G)$ has plenty of automorphisms, however, so $\cH(n)$ is not a scheme;
indeed, it is a {\it neutral gerbe}.
See  \cite{Smith}.

\begin{rema}[{\bf Landweber's criterion for flatness}]\label{LEFT} If $G:\Spec(R) \to \cM_\fg$
is any flat map, then there is a 2-periodic homology theory $E(R,G)$ with $E(R,G)_0 \cong R$
and with $G$ as the associated formal group. Landweber's Exact
Functor
Theorem gives an easily checked criterion to decide when a representable
morphism $\cN \to \cM_\fg$ is flat. 

Let $\cO_n$ be the structure sheaf of the substack $\cM(n)$ of $\cM_\fg$. Then the 
global section $v_n \in H^0(\cM(n),\omega^{\otimes p^n-1})$
defines an {\it injection} of sheaves
$$
\xymatrix{
0 \to \cO_n \ar[r]^-{v_n} \ar[r] &\omega^{\otimes p^n-1} }.
$$
This yields a short exact sequence
$$
\xymatrix{
0 \to \omega^{\otimes -(p^n-1)} \ar[r]^-{v_n} &\ar[r] \cO_n\ar[r] & j_\ast \cO_{n+1} \to 0.}
$$
This identifies $\omega^{\otimes -(p^n-1)}$ with the ideal defining 
the closed immersion of $\cM(n+1)$ in $\cM(n)$.

Now let $f:\cN \to \cM_\fg$ be a representable morphism of stacks and let
$$
\cN(n) = \cM(n) \times_{\cM_\fg} \cN \subseteq \cN.
$$
Then $\cN(n+1) \subseteq \cN(n)$ remains a closed immersion and,
if $f$ is flat, then
\begin{equation}\label{land-vn}
\xymatrix{
\cO_{\cN(n)} \ar[r]^-{v_n} & \omega^{\otimes p^n-1}
}
\end{equation}
remains an injection. Landweber's theorem now
says that this is sufficient; that is, if for all primes $p$ and
all $n$ the morphism of Equation (\ref{land-vn}) is an injection,
then the representable morphism $\cN \to \cM_\fg$ is flat. The original
source is \cite{Land}; in the form presented here, it appears in \cite{Hollander}
and \cite{Nau}.
\end{rema}

\begin{rema}[{\bf Chromatic stable homotopy theory}]\label{chrom}

The moduli stack $\cM_\fg$ of formal groups
has not been shown to be a derived stack and it may not be. One
technical difficulty is that we must  use the $fpqc$-topology on $\cM_\fg$ and there are
many flat maps. Nonetheless, the geometry of the stack $\cM_\fg$ has
been successfully used to predict theorems about large scale phenomena
in stable homotopy theory. Here are some of the basic results.

If $\cM_\fg$ could be lifted to a derived stack, there would be a descent
spectral sequence
$$
H^s(\cM_\fg,\omega^{\otimes t}) \Longrightarrow \pi_{2t-s}S.
$$
This spectral sequence exists: it is the Adams-Novikov Spectral Sequence. It remains our
most sensitive algebraic approximation to $\pi_\ast S$. The paper
\cite{MRW} initiated the modern era of calculations in stable homotopy theory;
much of the algebra there is driven by Morava's meditations on the geometry
of formal groups.

Another example is the Hopkins-Ravenel chromatic convergence theorem.
If we define 
$$
\cU(n) = \cM_\fg - \cM(n+1)
$$
to be the open substack classifying formal groups over schemes over $\ZZ_{(p)}$ of height
at most $n$, then we get an ascending chain of open substacks
$$
\QQ \otimes \cM_\fg \simeq \cU(0) \subseteq \cU(1) \subseteq \cU(2) \subseteq
\cdots \subseteq \ZZ_{(p)}\otimes \cM_\fg.
$$
This sequence is not exhaustive: the additive formal group of $\FF_p$ has infinite
height and does not give a point in any $\cU(n)$.

Now let $G_n:\Spec(R_n) \to \cM_\fg$ be a flat map classifying a formal group of exact height
$n$; this gives a 2-periodic homology theory $E(R_n,G_n)$. Let $L_n(-)$ be the localization
with respect to this homology theory. For example, $L_0X$ is rational localization
and $L_1X$ is localization with respect to $p$-local $K$-theory.
Then chromatic convergence \cite{RavOr} says that for  a  spectrum $X$
there is a tower
$$
\cdots \to L_2X \to L_1X \to L_0X$$
under $X$;  furthermore, if $\pi_k X = 0$ for $k$ sufficiently negative and $H_\ast (X,\ZZ_{(p)})$ is finitely generated
as a graded $\ZZ_{(p)}$-module, then
$$
X \longr \holim L_nX
$$
is  localization with respect to $H_\ast(-,\ZZ_{(p)})$.

Next we might like to decompose the $L_nX$. For this we use  the open inclusion
$i:\cU(n-1) \to \cU(n)$ and its closed complement $\cH(n) \subseteq \cU(n)$. Recall
$\cH(n)$ classifies formal groups of exact height $n$ and has a single geometric
point given by any formal group $\Ga_n$ over a field $\FF$ of characteristic $p$ of exact height $n$.
The classifying morphism $\Ga_n:\Spec(\FF) \to \cM_\fg$ is not flat, but the homology theory
$K(\FF,\Ga_n)$ exists nonetheless and is remarkably computable. These are
the {\it Morava $K$-theories}. If we write $L_{K(n)}(-)$ to denote localization
with respect to any of these theories (all the localizations are equivalent),
then there is a homotopy pull-back square
$$
\xymatrix{
L_nX \ar[r] \ar[d] & L_{K(n)}X \ar[d]\\
L_{n-1}X \ar[r] & L_{n-1}L_{K(n)}X.
}
$$
The square can be deduced from \cite{666}; the paper \cite{HovChrom}
contains a detailed analysis of a conjecture on how these squares behave.

We are thus left with calculating the pieces $L_{K(n)}X$. Here the theory is actually
fully realized. This is
the subject of the next subsection.
\end{rema}

\subsection{Deformations and the local Hopkins-Miller Theorem}\label{localHM}

We will need the language of deformation theory at several points; here it is used
to explain how one might compute $L_{K(n)}X$. Let $\cM$ be a stack and
$A_0/\FF$ be an $\cM$-object over a field $\FF$. Recall that an Artin local ring
$(R,\mm)$ is a local ring with nilpotent maximal ideal $\mm$. If $q:R \to \FF$
is a surjective morphism of rings, then a  {\it a deformation}
of $A_0$ to $R$ is an $\cM$-object $A$ and a
pull-back diagram
$$
\xymatrix{
A_0 \ar[r] \ar[d] & A \ar[d]\\
\Spec(\FF) \ar[r] & \Spec(R).
}
$$
Deformations form a groupoid-valued functor $\Def_\cM(\FF,A_0)$ on 
an appropriate category of Artin local rings.

If $\Ga$ is a formal group of finite height $n$ over a perfect field $\FF$, then Lubin-Tate
theory \cite{LT} says that the groupoid-valued functor $\Def_{\cM_\fg}(\FF,\Ga)$
is discrete; that is, the  natural
map
$$
\Def_{\cM_\fg}(\FF,\Ga) \to \pi_0\Def_{\cM_\fg}(\FF,\Ga)\defeq\hbox{Isomorphism classes in 
$\Def_{\cM_\fg}(\FF,\Ga)$}
$$
is an equivalence. Furthermore, $\pi_0\Def_{\cM_\fg}(\FF,\Ga)$ is
pro-represented by a complete local ring $R(\FF,\Ga)$; that is,
there is a natural isomorphism
$$
\pi_0\Def_{\cM_\fg}(\FF,\Ga) \cong \Spf(R(\FF,\Ga)).
$$
An appropriate choice of coordinate for the universal deformation
of $\Ga$ over $R(\FF,\Ga)$ defines an isomorphism
$$
W(\FF)[[u_1,\ldots,u_{n-1}]] \cong R(\FF,\Ga)
$$
where $W(-)$ is the Witt vector functor. If $\Aut(\FF,\Ga)$ is the group of
automorphisms of  the pair $(\FF,\Ga)$, then $\Aut(\FF,\Ga)$ acts on $\Def(\FF,\Ga)$
and hence on $R(\FF,\Ga)$. The formal spectrum $\Spf(R(\FF,\Ga))$ with the
action of $\Aut(\FF,\Ga)$ is called Lubin-Tate space; if we insist that $\FF$
is algebraically closed it is independent of the choice of the pair $(\FF,\Ga)$,
by Lazard's classification theorem. See \cite{GH} for the following result.

\begin{theo}[{\bf Local Hopkins-Miller}]\label{lhm}There is a $2$-periodic
commutative ring spectrum $E(\FF,\Ga)$ with $E(\FF,\Ga)_0 \cong R(\FF,\Ga)$
and associated formal group isomorphic to a universal deformation of $\Ga$.
Furthermore,
\begin{enumerate}

\item the space of all such commutative ring spectrum 
realizations of the universal deformation of $(\FF,\Ga)$ is
contractible; and

\item the group $\Aut(\FF,\Ga)$ acts on $E(\FF,\Ga)$ through maps of commutative
ring spectra.
\end{enumerate}
\end{theo}

The assignment $(\FF,\Ga) \mapsto E(\FF,\Ga)$ is actually a functor from a
category of height $n$ formal groups to commutative ring spectra.

In \cite{DH} Devinatz and Hopkins show that when $\FF$ contains enough
roots of unity, there is a weak equivalence
$$
L_{K(n)}X \simeq \holim_{\Aut(\FF,\Ga)} E(\FF,\Ga) \wedge X
$$
for all finite CW spectra and, in particular, for $X = S$.\footnote{There are technical
issues arising from the fact that $\Aut(\FF,\Ga)$ is a profinite group and that we must use
continuous cohomology. Again, see \cite{DH}.} This gives, for example,
a spectral sequence
$$
H^s(\Aut(\FF,\Ga),\omega_G^{\otimes t}) \Longrightarrow \pi_{2t-s}L_{K(n)}S
$$
where $G$ is a universal deformation of $\Ga$. If $p$ is large with respect to $n$,
this spectral sequence collapses and there are no extensions and the problem
becomes purely algebraic -- if not easy. For $n=1$ this happens if $p > 2$ and
if $n=2$, we need $p > 3$ and in these cases all the calculation have been done.
See \cite{Shi}. The case $p=2$ and $n=1$ is not hard; there has also been
extensive calculation at $p=3$ and $n=2$. In both cases $p$-torsion in $\Aut(\FF,\Ga)$
creates differentials. As a sample of the sort of large-scale periodic phenomena
we see, I offer the following classical result of Adams \cite{Ada1}. Let $p > 2$. Then $\pi_0L_{K(1)}S = \pi_{-1}L_{K(1)}S =
\ZZ_p$ and $\pi_kL_{K(1)}S$ is zero for all other $k$ unless
$$
k = 2p^ts(p-1) -1,\qquad (s,p) = 1
$$
and then
$$
\pi_kL_{K(1)}S = \ZZ/p^{t+1}.
$$
All the elements in positive degree come from $\pi_\ast S$ itself; in fact,  the natural
map $\pi_\ast S \to \pi_\ast L_{K(1)}S$ splits off the image
of the classical $J$-homomorphism.

\section{Construction and deconstruction on $\mell$.}

\subsection{A realization problem}

If $X \to \mell$ classifies a generalized elliptic curve $C \to X$, let 
$\cI(e) \subseteq \cO_C$ be the ideal sheaf defining the identity section $e:X \to C$.
Then $\cI(e)/\cI(e)^2$ is an $e_\ast \cO_X$-module sheaf and thus determines
a unique $\cO_X$-module sheaf $\omega_C$ on $X$. This sheaf is locally free
of rank 1. Even if $C$ is singular, it is never singular at $e$ and the smooth
locus of $C$ has the structure of an abelian group scheme of dimension 1 over
$X$; hence, $\omega_C$ is isomorphic to the invariant differentials of $C$.
The assignment
$$
\omega_\Ell(C:X \longr \mell) = \omega_C
$$
defines an invertible sheaf on $\mell$. 

More is true. If $C_e$ denotes the formal completion of $C$ at $e$, then $C_e$
is a smooth 1-parameter formal Lie group and $\omega_C \cong \omega_{C_e}$.
This data defines a morphism
$$
q:\mell \longr \cM_\fg
$$
with the property that $q^\ast \omega = \omega_\Ell$; for this reason I will simply
write $\omega_\Ell$ as $\omega$, leaving the base stack to determine which sheaf I mean.

If $C$ is a generalized elliptic curve over an affine scheme $X$ and $\omega_C$ has a trivializing
global section $\eta:\cO_X \to \omega_C$, then there is a choice of closed immersion
$C \to \PP^2$ over $X$ with $e$ sent to $[0,1,0]$. This closed immersion is defined by a Weierstrass equation
$$
Y^2Z + a_1XYZ + a_3YZ^2 = X^3 + a_2X^2Z + a_4 XZ^2 + a_6 Z^3
$$
which is normally written as
$$
y^2 + a_1xy + a_3y = x^3 + a_2x^2 + a_4 x + a_6
$$
since there is only one point where $Z=0$. This closed immersion is not unique, but is
adapted to $\eta$ in the sense that $z = -x/y$ reduces to $\eta$ in $\omega_C$.
Other choices of $\eta$ and the associated adapted immersion produce a different
immersion; all such immersions differ by a projective transformation.
Given $\eta$ we see that the invariant differentials for the associated formal
group $C_e$ form the trivial sheaf; hence $C_e$ can be given a coordinate. Once we have chosen
a coordinate $t$ and an immersion adapted to $\eta$ we have an equation
$$
z= t + e_1t^2 + e_2t^3 + \cdots
$$
since $z$ is also a coordinate. Then 
we can uniquely specify the adapted
immersion by requiring $e_i=0$ for $i \leq 3$. From this
we can conclude that the morphism $\mell \to \cM_\fg$ is representable. This argument
is due to Hopkins.

This morphism is also flat, by Landweber's criterion \ref{LEFT}. Indeed, in this setting the
global section
$$
v_1 \in H^0(\FF_p \otimes \mell,\omega^{\otimes p-1})
$$
is the {\bf Hasse invariant}. A generalized elliptic curve $C$ over $X$ over $\FF_p$
is {\bf ordinary} if $v_1(C):\cO_X \to \omega_C$ is an isomorphism; on the other hand,
if $v_1(C) = 0$,
then $C$ is {\bf supersingular} and automatically smooth.
We define $\cM_\ss \subseteq \FF_p \otimes \mell$ to be the
closed substack defined by the vanishing of $v_1$ and then, over $\cM_\ss$, 
global section 
$$
v_2:\cO_\ss \longr \omega^{\otimes p^2-1}
$$
is an isomorphism. This is a rephrasing of the statement that formal group of a
supersingular curve has exact height $2$.

\begin{rema}[{\bf The realization problem}]\label{real-prob}
We can now ask the following question. Suppose $\cN \to \cM_\fg$ is a representable
and flat morphism from a Deligne-Mumford stack to the moduli stack of formal groups.
A realization of $\cN$ is a derived Deligne-Mumford stack $(\cN,\cO)$ with
\[
\pi_{k}\cO \cong \brackets{\omega_\cN^{\otimes t},}{k=2t;}{0,}{k=2t+1.}
\]
Does $\cN$ have a realization? If so, how many are there? Better, what's the homotopy
type of the space of all realizations?

This question is na\"\i ve for a variety of reasons. A simple one is that if $\cO$ exists,
then $\pi_\ast \cO$ will be a sheaf of graded commutative rings with power operations
and, in general, there is not enough data in formal groups to specify those operations.
In the case of the local Hopkins-Miller theorem \ref{lhm} these operations are determined
by the subgroup structure of the formal groups. See \cite{AHS2} and \cite{Rez2}. In
the case of $\mell$, they will be determined by the Serre-Tate theorem and the
subgroup structure of elliptic curves. 

Work
of Mark Behrens and Tyler Lawson \cite{BL} solves the realization problem for
certain Shimura varieties, which are moduli stacks of highly structured
abelian varieties. The extra structure is needed to get formal groups
of higher heights.
\end{rema}

The realization problem is essentially a $p$-adic question. If $X$ is a scheme
or a stack, let $j_n:X(p^n) \to X$ be the closed immersion defined by the vanishing
of $p^n$.  The formal completion of $X$ at $p$ is the colimit sheaf $X^\wedge_p = \colim X(p^n)$;
as a functor, $X^\wedge_p$ is the restriction of $X$ to rings in which $p$ is nilpotent; thus
$X^\wedge_p$ is a formal scheme over $\Spec(\ZZ)^\wedge_p = \Spf(\ZZ_p)$. 
If $\cF$ is a derived module sheaf on $X$, the derived completion of $\cF$ is
$$
\cF_p = \holim (j_n)_\ast j_n^\ast \cF.
$$
Then if $\cF$ is  a
derived sheaf of modules on $\cN$ then there is a homotopy pull-back square
\begin{equation}\label{arithmetic}
\xymatrix{
\cF \ar[r] \ar[d] & \prod_p \cF_p\ar[d]\\
\QQ \otimes \cF \ar[r] & \QQ \otimes (\prod_p \cF_p)
}
\end{equation}
where $\QQ \otimes (-)$ is rational localization. Here I am regarding $\cF_p$ as a sheaf
on $X$; however, $\cF_p$ is determined by its restriction to $X_p$.

The category
of commutative ring spectra over the rational numbers is Quillen equivalent to the category
of differential graded algebras over $\QQ$; hence, many questions in derived algebraic geometry
over $\QQ$ become classical.

\begin{exem}[{\bf $\mell$ over $\QQ$}]\label{mell-at-q}
If $S_\ast$ is a graded ring with $S_0 = R$, then the grading 
defines an action of the multiplicative group $\GG_m$ on $Y = \Spec(S_\ast)$. Define
$\Proj(S_\ast)$ to be the quotient stack over $R$ of this action on $Y - \{0\}$ where $\{0\} \subseteq Y$
is the closed subscheme defined by the ideal of elements of positive degree. Then
$$
\QQ \otimes \mell \cong \Proj(\QQ[c_4,c_6])
$$
and the graded sheaf of rings $\QQ[\omega^{\pm 1}]$ is a dga with trivial
differential. This defines the derived scheme $\QQ \otimes \mell$ and the sheaf
$\QQ \otimes \cO$. Similarly $\QQ \otimes (\mell)_p = \Proj(\QQ_p[c_4,c_6])$
and the map across the bottom in Diagram (\ref{arithmetic}) will be the obvious one.
I note that we can define $\ZZ[1/6] \otimes \mell$ as $\Proj(\ZZ[1/6][c_4,c_6])$; however,
we can't define the derived scheme this way as the homotopy theory of commutative
ring spectra over $\ZZ[1/6]$ is not equivalent to a category of dgas.
\end{exem}

\subsection{Lurie's theorem and $p$-divisible groups}

In this section, we pick a prime $p$ and work
over $\Spf(\ZZ_p)$; that is, $p$ is implicitly nilpotent in all our rings. I will
leave this out of the notation.

\begin{defi}\label{pdiv} Let $R$ be a ring and $G$ a sheaf of abelian groups
on $R$-algebras. Then $G$ is a {\bf $p$-divisible group} of {\bf height $n$}
if 
\begin{enumerate}

\item $p^k:G \to G$ is surjective for all $k$;

\item $G(p^k) = \mathrm{Ker}(p^k:G \to G)$ is a finite and flat group
scheme over $R$ of rank $p^{kn}$;

\item $\colim G(p^k)\cong G$.
\end{enumerate}
\end{defi}

\def\et{{{\mathrm{{et}}}}}
\def\for{{{\mathrm{for}}}}

\begin{rema}\label{p-div-rem} 1.) If $G$ is a $p$-divisible group, then
completion at $e \in G$ gives an abelian formal group $G_\for \subseteq G$,
not necessarily of dimension $1$. The quotient $G/G_\for$ is 
\'etale over $R$; thus we get a natural short exact sequence
$$
0 \to G_\for \to G \to G_\et \to 0.
$$
This is split over fields, but not in general.

2.) If $C$ is a smooth elliptic curve, then $C(p^\infty) = \colim\ C(p^n)$ is
$p$-divisible of height $2$ with formal part of dimension 1.

3.) If $G$ is a $p$-divisible group over a scheme $X$, the function
which assigns to each geometric point $x$ of $X$ the height of the
fiber $G_x$ of $G$ at $x$ is constant. This is not true of formal groups, as
the example of elliptic curves shows. Indeed, if $G$ is $p$-divisible
of height $n$ with $G_\for$ of dimension 1, then the height of $G_\for$
can be any integer between $1$ and $n$.
\end{rema}

\begin{defi}\label{mod-p-div} Let $\cM_p(n)$ be the moduli stack
of $p$-divisible groups of
height $n$ and with $\mathrm{dim}\ G_\for = 1$.
\end{defi}

\begin{rema}\label{p-div-rems}The stack $\cM_p(n)$ is not an algebraic stack,
but rather pro-algebraic in the same way $\cM_\fg$ is pro-algebraic.
This can be deduced from the material in the first chapter of \cite{Messing}.
\end{rema}

There is a morphism of stacks $\cM_p(n) \to \cM_\fg$ sending
$G$ to $G_\for$. This morphism is {\it not} representable.
By definition, there is a factoring of  this map as
$$
\cM_p(n) \longr \cU(n) \longr \modfg
$$
through the open substack of formal groups of height at most $n$.

We now can state Lurie's realization result \cite{Lurell}.
Since we are working over $\ZZ_p$, one must take care with 
the hypotheses: the notions of algebraic stack
and \'etale must be the appropriate notions
over $\Spf(\ZZ_p)$. 

\begin{theo}[{\bf Lurie}]\label{mainthm} Let $\cM$ be an Deligne-Mumford stack
equipped with a formally \'etale morphism
$$
\cM \longr \cM_p(n).
$$
Then the realization problem for the composition
$$
\cM \longr \cM_p(n) \longr \cM_\fg
$$
has a canonical solution; that is, the space of all solutions 
has a preferred basepoint.
\end{theo}

It is worth emphasizing that this theorem uses the local Hopkins-Miller theorem
\ref{localHM} in an essential way.

\begin{exem}[{\bf Serre-Tate theory}]\label{ST} As addendum to this theorem,
Lurie shows that the morphism $\epsilon:\cM \to \cM_p(n)$ is formally \'etale if
it satisfies the Serre-Tate theorem. This applies to the open
substack $\cM_{1,1}$ of $\mell$ and we recover the main theorem
\ref{mn-thm}, at least for smooth elliptic curves.  Then, in \cite{Lurell}, Lurie
gives an argument extending the result to the compactification
$\mell$.

The Serre-Tate theorem asserts an equivalence of deformation groupoids;
compare the Lubin-Tate result of Section \ref{localHM}.
Let $\cM$ be a stack over $\cM_p(n)$ and $A_0/\FF$ be an $\cM$-object
over a field $\FF$, necessarily of characteristic $p$ since we are working
over $\Spf(\ZZ_p)$. Then the deformations form a groupoid functor $\Def_\cM(\FF,A_0)$ on 
Artin local rings. The Serre-Tate theorem holds if the evident morphism
$$
\Def_\cM(\FF,A_0) \longr \Def_{\cM_p(n)}(\FF,\epsilon A_0)
$$
is a equivalence. This result holds for elliptic curves, but actually
in much wider generality. See \cite{Messing}.
\end{exem}

\begin{rema}[{\bf Deformations of $p$-divisible groups}]\label{Def} 
We discussed the deformation theory of formal groups and Lubin-Tate
theory in Section \ref{localHM}; the theory is very similar for
$p$-divisible groups. Let $G$ be a 
$p$-divisible group over an algebraically closed field $\FF$. Then we have 
split short exact sequence 
$$
0 \to G_\for \to G \to G_\et \to 0.
$$
Since $G_\et$ has a unique deformation up to isomorphism, by the definition
of \'etale, the deformations of $G$ are determined by the deformations of
$G_\for$ and an extension class. From this
it follows that the groupoid-valued functor $\Def_{\cM_p(n)}(\FF,G)$
is discrete and $\pi_0\Def_{\cM_p(n)}(\FF,G)$ is pro-represented by
$$
R(\FF,G_\for)[[t_1,\cdots,t_{n-h}]] \cong W(\FF)[[u_1,\cdots,u_{h-1},t_1,\cdots,t_{n-h}]].
$$
Note that this is always a power series in $n-1$ variables.
Using this remark it is possible to give a local criterion for when a morphism
of stacks $\cM \to \cM_p(n)$ is \'etale. It is in this guise that Lurie's theorem
appears in \cite{BL}.
\end{rema}

\subsection{Decomposing the structure sheaf}

Both the original Hopkins-Miller argument and Lurie's derived
algebraic geometry construction of the derived $\mell$
rely on a decomposition of the moduli stack
of elliptic curve into its ordinary and supersingular components. 
Specifically, let $\cM_\ss \subseteq \mell$ be the closed substack of
supersingular curves defined by the vanishing of both a fixed prime $p$ and the
Hasse invariant. There are inclusion of closed substacks
$$
\cM_\ss \subseteq \FF_p \otimes \mell \subseteq \mell;
$$
from this we get an inclusion of formal substacks $\hatss  \subseteq (\mell)^\wedge_p$.
Let $\hatord \subseteq  (\mell)^\wedge_p$ be the open complement of $\hatss$; thus,
$\hatord$ classifies ordinary elliptic curves over rings $R$ in which $p$ is nilpotent.
If $X$ is a sheaf on $\mell$, we obtain completions of $X$
from each of these formal stacks, which we write $L_pX$, $L_\ss X$, and $L_\ord X$,
respectively.

\begin{theo}
\label{2nd-assemby}
There is a homotopy pull-back diagram
of sheaves of $E_\infty$-ring spectra on $\mell$
$$
\xymatrix{
L_p \cO \ar[r] \ar[d] & L_\ss \cO \ar[d]\\
L_\ord \cO \ar[r] & L_\ord L_\ss \cO.
}
$$
\end{theo}

This statement belies history: the original Hopkins-Miller proof used obstruction theory
to build the sheaves $L_\ord \cO$ and $L_\ss \cO$ and the map across the bottom
of this diagram, thus building $L_p \cO$. Lurie, on the other hand, produces a candidate
for $L_p\cO$, and then must show it has the right homotopy type. This is done by
analyzing the pieces in this square separately. In either case, the pieces $L_\ord \cO$
and $L_\ss \cO$ have intrinsic interest and the calculations of their homotopy types 
calls on classical calculations with modular forms.

\subsection{The supersingular locus}

The description of $L_\ss\cO$ begins with an alternative description of the stack
$\mss$ of supersingular curves. This is standard theory for elliptic curves; see
\cite{KM}, Chapter 12. Since the vanishing of the Hasse invariant of $C$ 
depends only on the isomorphism class of the curve $C$, it depends only on the $j$-invariant
of $C$ as well. Write
$$
j:\FF_p \otimes \mell \longr \FF_p \otimes \PP^1=\PP^1
$$
with $j(C) = [c_4^3(C),\Delta(C)]$ and call the point $[1,0] \in \PP^1$  the point at infinity.
If $\Delta(C)$ is invertible -- that is, $C$ is an elliptic curve -- we also write
$$
j(C) = c_4^3(C)/\Delta(C) \in \AA^1_j = \PP^1 - \{\infty\}.
$$
A point in $\PP^1$ will be ordinary if it is the image of an ordinary curve; supersingular
otherwise.

\begin{theo}\label{ss-points-in-p1} Fix a prime $p$. The point at $\infty$ in $\PP^1$ is ordinary
and there is a separable polynomial $\Phi(j) \in \FF_p[j]$ so that $C$ is supersingular
if and only if $j(C)$ is a root of $\Phi$. All the roots of $\Phi$ lie in $\FF_{p^2}$.
\end{theo}

For example, if $p=2$ or $3$, then $\Phi(j) = j$.
The degree of $\Phi$ is $[(p-1)/12] + \epsilon_p$ where $\epsilon_p$ is $0$, $1$, or $2$
depending on the prime.

\begin{prop}\label{mss-described} Let $R_\ss = \FF_p[j]/(\Phi(j))$.
There is a
$2$-category pull-back 
\[
\xymatrix{
\mss \ar[r] \ar[d] & \FF_p \otimes \mell\ar[d]^j\\
\Spec(R_\ss) \ar[r] & \PP^1.
}
\]
Furthermore, there is a supersingular curve $C_\ss$ over $\FF_{p^2} \otimes R_\ss$ which
gives an equivalence of stacks
$$
B\Aut(C_\ss) \simeq \FF_{p^2} \otimes \mss.
$$
\end{prop}

Here I have written $\Aut(C_\ss)$ for $\Aut(C_\ss/(\FF_{p^2} \otimes R_\ss))$.
If $G$ is a group scheme, $BG$ is the moduli stack of $G$-torsors.

To be concrete,  write $C_\ss = \coprod C_a$ where $C_a$ is
a representative for the isomorphism class of supersingular curves over
$\FF_{p^2}$ with $j(C_a)=a$.
Then
\begin{equation}\label{equiv-decomp}
\FF_{p^2} \otimes \mss \simeq \coprod B\Aut(C_a/\FF_{p^2}).
\end{equation}

The algebraic groups $\Aut(C_a/\FF_{p^2})$ are all known. See \cite{Silver}, among
many other sources.
The difficulties and the interest lie at the primes $2$ and $3$, 
where the $\FF_{p^2}$ points
of $\Aut(C_0/\FF_{p})$ have elements of order $p$. For example, if $p=3$,
$$
\Aut(C_0/\FF_9) \cong \mu_4 \rtimes \ZZ/3.
$$
where $\mu_4$ is the 4th roots of unity.

Let $C:\Spec(B) \to \mell$ be \'etale and $I = I(B,C) \subseteq B$ the ideal generated
by $p$ and the Hasse invariant. Then $B/I$ is a separable $\FF_p$-algebra
and $q^\ast C_e$ is a formal group of exact height $2$. The evident extension of the local
Hopkins-Miller theorem \ref{lhm} to separable algebras gives a sheaf $\cE_\ss$ of commutative ring
spectra with
$$
\cE_\ss(C:\Spec(B) \to \mell) \simeq E(B/I,q^\ast C_e)
$$
and we have $L_\ss \cO \simeq \cE_\ss$.

It is straightforward from here to understand the homotopy type of the global sections
of $L_\ss \cO$; indeed, we have
$$
R\Ga L_\ss \cO \simeq \holim_G E(\FF_{p^2} \otimes R_\ss, (C_\ss)_e)
$$
where $G = \Gal(\FF_{p^2}/\FF_p) \rtimes \Aut(C_\ss)$.

The homotopy groups of this spectrum have been computed.
If $p > 5$ this is fairly easy. If $p=3$ the most explicit source is \cite{GHMR},
and if $p=2$ it is implicit in \cite{HopMah} although that source needs to be combined with \cite{Bau}
to get complete answers.

\subsection{The ordinary locus}

For the sheaf $L_\ord \cO$ we use the map $\mell \to \cM_\fg$ classifying the associated
formal group  to make preliminary computations. By construction, this morphism
restricts to a morphism $\hatord \to \cU(1)$ to the open substack of formal groups
of height 1. The map $\cU(1) \to \Spf(\ZZ_p)$ is formally \'etale and has section
$g:\Spf(\ZZ_p) \to \cU(1)$
classifying the multiplicative formal group $\hatgm$. The map $g$ is pro-Galois with 
Galois group $\Aut(\hatgm) = \ZZ_p^\times$, the units in the $p$-adics.  We use
this cover not to define $L_\ord \cO$ directly, but to first specify the resulting
homology sheaf $K_\ast \cO$, where $K_\ast$ is $p$-complete complex
$K$-theory.\footnote{Because we are working
$p$-adically, $K_\ast X = \pi_\ast L_{K(1)}K \wedge X$, a completion of the usual
$K$-theory. The issues here are very technical, but carefully worked out
in \cite{666}.}
From this we then can construct $L_\ord \cO$.

Let $R$ be a $p$-complete ring and let 
$C:\Spec(R) \to \mell$ be an \'etale morphism. Consider the following
diagram, where both squares are 2-category pull-backs
$$
\xymatrix{
\Spf K_0\cO(R,C) \ar[d]\ar[r] & \Spf(V) \ar[d]\ar[r] & \Spf(\ZZ_p) \ar[d]^{\hatgm}\\
\Spec(R) \ar[r] & \mell \ar[r] & \cM_\fg.
}
$$
This defines $K_0\cO(R,C)$ and the ring $V$. The latter ring is Katz's ring
of {\it divided congruences} \cite{Kat} -- and, as it turns out, $K_0\tmf$. Note that $\Spf(V)$
solves the moduli problem which assigns to each ring $A$ with $p$ nilpotent
the set of pairs $(C,\phi)$ where $C$ is an elliptic curve over $A$ and $\phi:\hatgm \to C_e$
is an isomorphism. By construction $V \to K_0\cO(R,C)$ is \'etale. We extend
the resulting sheaf of rings $K_0\cO$ to a graded sheaf $K_\ast \cO$ by
twisting by the various powers of the sheaf of
invariant differentials, as in Equation (\ref{inv-diff-defn-ext}).

The sheaf of rings $K_0\cO$ has a great deal of structure. First, it
comes equipped with an action of $\ZZ_p^\times$ from the automorphisms
of $\hatgm$. This gives the action of the Adams operations in $K$-theory. 
Second, and more subtly, there is an extra ring operation $\psi$
which is a lift of the Frobenius and commutes with the action of $\Aut(\GG_m)$. 
To construct $\psi$, we use that if $C$ is an ordinary elliptic curve, then
the kernel of $p:C_e \to C_e$ defines a canonical subgroup $C(p) \to C$
of order $p$.
Then $\psi:V \to V$ is defined by specifying the natural transformation
on the functor $V$ 
represents:
$$
\psi(C,\phi:\hatgm \to C) = (C/C(p),\hatgm \cong \hatgm/\hatgm(p) \mathop{\to}^{\bar{\phi}} (C/C(p))_e).
$$
If $C$ is an ordinary curve over an $\FF_p$-algebra $A$, then $C/C(p) = C^{(p)}$; this 
explains why
$\psi:V \to V$ is a lift of the Frobenius. To extend $\psi$ to all of
$K_0\cO$, note that 
since $V \to K_0\cO(R,C)$ is \'etale there is a unique morphism $\psi:K_0\cO(R,C)
\to K_0\cO(R,C)$ lifting the Frobenius and extending $\psi$ on $V$.

Since $K_0\cO$ is torsion-free $\psi(x) = x^p + p\theta(x)$ for some unique
operator $\theta$. All this structure extends in a unique way to $K_\ast \cO$,
making this a sheaf of {\it $theta$-algebras}. By work of McClure \cite{McL}, this is exactly
the algebraic structure supported by $K_\ast X$ when $X$ is a commutative ring
spectrum. The existence of $L_\ord \cO$ is now guaranteed by the following result.

\begin{theo}\label{lk1-tmf} The space of all sheaves of $L_{K(1)}$-local ring spectra
$X$ with $K_\ast X \cong K_\ast \cO$ as sheaves of theta-algebras is non-empty
and connected.
\end{theo}

Any of these sheaves -- they are all homotopy equivalent -- gives a model for
$L_\ord \cO$. Lurie's work does better: it gives a canonical model.

Originally, Theorem \ref{lk1-tmf} was proved by an obstruction theory argument; the obstruction groups
are computed using an appropriate cotangent complex for the sheaf $K_0\cO$. This
argument uses in an essential way that $\mell$ is smooth of dimension 1 over $\ZZ$.
At the crucial primes $2$ and $3$ there is a very elegant construction originating
with Hopkins of $L_{K(1)}\tmf$ itself which short-circuits the obstruction theory.
See the paper by Laures \cite{Lau}.

\begin{rema}\label{homotopy-of-ord-tmf}The homotopy groups of the derived global
sections $R\Gamma(L_\ord\cO) \simeq L_{K(1)}\tmf$ have been computed. For example,
at the prime $2$, we have
$$
\pi_\ast R\Gamma(L_\ord\cO) \cong (\ZZ_2[1/j])^{\wedge}_2[\eta,v,b^{\pm 1}]/I
$$
where the comleted polynomial ring on the inverse of the $j$-invariant
is in degree $0$ and $\eta$, $v$, and $b$ have
degrees $1$, $4$, and $8$ respectively. The ideal $I$ is generated by the relations
\begin{align*}
2\eta = \eta^3 &= v\eta = 0\\
v^2 &= 2b.
\end{align*}
There is a map $R\Gamma(L_\ord\cO) \to KO$ to the spectrum of $2$-completed real $K$-
theory which, in homotopy, is the quotient by the ideal generated by $1/j$. 
\end{rema}

\begin{rema}\label{finish} To complete the construction of $L_p\cO$, we must produce the map
$L_\ord \cO \to L_\ord L_\ss \cO$. To do this, we calculate $K_\ast L_\ss\cO$ as
a theta-algebra and again use obstruction theory. The complication
is that the lift $\psi$ of the Frobenius on $K_\ast L_\ss\cO$ is determined
by the $E_\infty$-ring structure on $L_\ss \cO$. For the algebra to work, 
we need to check that this is the same lift as that determined by the subgroup structure
of an appropriate elliptic curve. To obtain the elliptic curve, we use the Serre-Tate
theorem to show that the universal deformation of the formal
group of a supersingular curve is actually the $p$-divisible group of an elliptic curve.  
Then we can apply  the general theory of power operations in complex
orientable homology theories as developed by Ando \cite{AHS2}
and Rezk \cite{Rez2}. I have to thank Charles Rezk for explaining this to me.
Lurie's construction avoids this question, because the map already exists.
\end{rema}

\subsection{Derived modular forms and duality}

The descent spectral seqeunce
\begin{equation}\label{descentss}
H^s (\mell,\omega^{\otimes t}) \Longrightarrow \pi_{2t-s}\tmf
\end{equation}
has been completely calculated and the homotopy groups of $\tmf$ exhibit a very
strong form of duality not present in the cohomology groups of $\mell$. Let me try to give
some flavor of the results. The first observation is that, by \cite{Del},
we have an isomorphism
$$
M_\ast = \ZZ[c_4,c_6,\Delta] \cong H^0(\mell,\omega^{\otimes \ast})
$$
from the ring of modular forms of level 1 to the global sections.

Next, the stack $\mell$ is smooth of dimension $1$ over $\ZZ$ and and the
cotangent sheaf is identified by the isomorphism
$$
\Omega_{\mell/\ZZ} \cong \omega^{\otimes -10}.
$$
I learned this from Hopkins, and it can be deduced from \cite{KM} Chapter 10.
There is an isomorphism $H^1(\mell,\omega^{\otimes -10})\cong \ZZ$ and
a multiplication homomorphism
\begin{equation}\label{pair6}
H^1(\mell,\omega^{\otimes -k-10}) \otimes H^0(\mell,\omega^{\otimes k})
\to 
H^0(\mell,\omega^{\otimes -10}) = \ZZ
\end{equation}
If we invert the primes $2$ and $3$, $H^s(\mell,\omega^{\otimes \ast})=0$
for $s >1$ and the multiplication of Equation (\ref{pair6})
becomes a perfect pairing. Thus, if we invert $6$ the spectral sequence of
(\ref{descentss}) collapses,
there can be no possible extensions, and both the coherent cohomology of
$\mell$ and the homotopy groups of $\tmf$
exhibit a strong form of Serre-type duality. This is not very surprising as, when $6$ is inverted,
$\mell$ is isomorphic to the projective stack obtained from the graded ring $M_\ast$.
See Remark \ref{mell-at-q}. 

Over the integers, however, the behavior is more subtle. The bookkeeping at the prime
$2$ is difficult, so let me say what happens when $2$ is inverted. The presence of
an element of order $3$ in the automorphisms of the supersingular curve at the
prime $3$ forces higher cohomology. There is an injection
\begin{equation}\label{doublem}
\grm \defeq M_\ast [\alpha,\beta]/R \longr H^\ast (\mell,\omega^{\otimes \ast})
\end{equation}
where $\alpha \in H^1(\mell,\omega^{\otimes 2})$ and $\beta \in H^2(\mell,\omega^{\otimes 6})$
and $R$ is the ideal of relations
$$
3\alpha = 3\beta =\alpha^2=c_i\alpha = c_i\beta=0.
$$
The element $\alpha$ is the image of the Hasse invariant under the boundary map
$$
H^0(\FF_3 \otimes \mell,\omega^{\otimes 2}) \to H^1(\mell,\omega^{\otimes 2})
$$
and $\beta$ is the Massey product $\langle \alpha,\alpha,\alpha\rangle$.
Both $\alpha$ and $\beta$ arise from the homotopy groups of spheres. The ring
$\grm$ is actually the coherent cohomology of a moduli stack of Weierstrass curves.
See \cite{Bau}.

We now define a bigraded $\grm$-module $\otherh$ by the short exact sequence
of $\grm$-modules
\begin{equation}\label{ses-in-coh}
0 \to \grm \to H^\ast(\mell,\omega^{\otimes \ast}) \to \otherh \to 0.
\end{equation}
We then need to write down the $\grm$-module $\otherh$. It will turn out
that in any given bidegree one of (and possibly both)
$\grm$ and $\otherh$ is zero, so we get an unamibigous splitting of bigraded groups
$$
H^\ast(\mell,\omega^{\otimes \ast}) = \grm \oplus \otherh.
$$
It is not quite split as $\grm$-modules as mulitplication by $\Delta$ will link the two pieces.

To describe the modules $\otherh$, we first note that
the pairing of (\ref{pair6}) is no longer a perfect pairing with only $2$ inverted: it only induces an
injection. First one computes $C^1_t = H^1(\mell,\omega^{\otimes t})$ if $t < 0$ and $C^1_t=0$ otherwise.
Now let $K^1_\ast \subseteq C^1_\ast$ be the kernel of multiplication by $\beta$. Then,
the pairing induces a morphism of short exact sequences of graded
$M_\ast$-modules
$$
\xymatrix@C=15pt{
\ \ \ \ \ \ \ 0 \ar[r]& K^1_\ast \ar[d]_\cong \ar[r] & C^1_\ast \ar[d]\ar[r] &
\Hom_{\FF_3}(\FF_3[\Delta],\FF_3)\ar[d]\ar[r]&  0\\
0 \ar[r]& \Hom_{\ZZ}(M_\ast,3\ZZ[1/2]) \ar[r] & \Hom_{\ZZ}(M_\ast,\ZZ[1/2])\ar[r]&
\Hom_{\FF_3}(M_\ast/3,\FF_3) \ar[r] & 0
}
$$
where the right vertical inclusion is induced by the quotient by the ideal $(c_4,c_6)$ of $M_\ast$.

To complete the description of $\otherh$ 
we extend the top short exact sequence of our diagram
of $M^{\ast}$-modules to a short exact sequence of $M^{\ast}_{\ast}$-modules by
$$
0 \to K^1_\ast \to \otherh \to \grm \otimes_{M_\ast} \Hom_{\FF_3}(\FF_3[\Delta],\FF_3) \to 0.
$$
The extension of $\grm$-modules in the exact sequence (\ref{ses-in-coh}) is now determined
by the requirement that multiplication by $\Delta$ is surjective in positive cohomological
degrees. In these degrees multiplication by $c_4$ and $c_6$ is necessarily zero,
for degree reasons.

There must be differentials in the descent spectral sequence (\ref{descentss}). A classical
relation, due to Toda, states that $\alpha\beta^2=0$ in homotopy; this forces
$d_5(\Delta) = \pm \alpha\beta^2$ in the spectral sequence. This differential and elementary
considerations dictates the entire spectral sequence. Let $DM_\ast \subset M_\ast$
be the subring generated by $c_4,c_6, 3\Delta,3\Delta^2$, and $\Delta^3$. The
inclusion $DM_\ast \subseteq \tmf_\ast$ extends to a split inclusion
$$
DM^{\ast}_{\ast} = DM_\ast[\alpha,\beta,x]/DR \subseteq \tmf_\ast
$$
where $DR$ is the ideal of relations
\begin{align*}
3\alpha = 3\beta = 3x &= c_i\alpha = c_i\beta = c_ix=0\\
\alpha^2 = \beta^5&= \alpha\beta^2=x\beta^2=0\\
\alpha x &= \beta^3 
\end{align*}
The class $x$ is the Toda bracket $\langle \alpha,\alpha,\beta^2\rangle$ and is detected
by $\alpha\Delta$ in the spectral sequence. The inclusion $DM^{\ast}_{\ast}
\subseteq \tmf_\ast$ of $DM^{\ast}_{\ast}$-modules has cokernel $D\otherh$.
There is a short exact
sequence of $DM^{\ast,\ast}$-modules 
$$
0 \to K^1_\ast \to D\otherh \to D\grm \otimes_{DM_\ast} \Hom_{\FF_3}(\FF_3[\Delta^3],\FF_3) \to 0.
$$
Note that there is a quotient $DM_\ast \to \FF_3[\Delta^3]$; this defines the module
structure needed for the tensor product.

\begin{rema}[{\bf Duality for $\tmf$}]\label{duality}At this point a new duality emerges --
one not apparent before completing these homotopy theory calculations. 
A first remarkable feature is that $D\otherh$ is in degrees less than $-20$ and
$D\grm$ is concentrated in non-negative degrees. In particular $D\grm$ is the homotopy
groups of the connected cover of $\tmf$.

A second feature is that $D\grm$ and $D\otherh$ are almost dual as $D\grm$-modules.
There are a number of ways to make this precise; a simple one is
to say that for all primes $p$ there is a homomorphism $\tmf_{-21} \to \FF_p$ so the
induced map given by ring multiplication
$$
\pi_k(\tmf/p) \to \Hom(\pi_{-k-21}\tmf/p,\FF_p)
$$
is an isomorphism. Here $\tmf/p$ is derived global sections of the topological structure
sheaf on $\FF_p \otimes\mell$.

This duality has an elegant homotopy theoretic explanation, given by Mahowald
and Rezk \cite{MR1}; in that source,
there is also a simple explanation of the number $-21$. In derived geometry, this
number surely appears because the dualizing sheaf is $\omega^{\otimes -10}$. However,
I know of no explanation for the duality from the 
point of view of derived algebraic geometry.
\end{rema}

\end{document}